\documentclass[a4paper,11pt,twoside]{article}
\usepackage{amssymb}
\usepackage{pstricks}
\usepackage{pst-plot}
\usepackage{pst-node}
\usepackage{fancyhdr}

\newtheorem{theorem}{Theorem}
\newtheorem{lemma}{Lemma}
\newtheorem{proposition}{Proposition}

\newtheorem{corollary}{Corollary}
\newtheorem{remark}{Remark}

\title{Spectral decomposition and Gelfand's theorem}
\author{A. Driouich, O. El-Mennaoui\\ Agadir university\\ and\\ M. Jazar\footnote{
Corresponding author.\newline The third author is
supported by a grant from the Lebanese University.}\\
Lebanese university}
\begin{document}

\maketitle

\begin{center}
{\bf\small Abstract}

\vspace{3mm} \hspace{.05in}\parbox{4.5in} {{\small In this paper
we are interested in spectral decomposition of an unbounded
operator with discrete spectrum. We show that if $A$ generates a
polynomially bounded $n$-times integrated group whose spectrum set
$\sigma(A)=\{i\lambda_k;\,k\in\mathbb{Z}^* \}$ is discrete and
satisfies $\sum \frac{1}{|\lambda_k|^\ell\delta_k^n}<\infty$ ($n$
and $\ell$ nonnegative integers), then there exists projectors
$(P_k)_{k\in\mathbb{Z}^*}$ such that $\sum P_kx=x$ ($ x\in D(
A^{n+\ell})$), where $\delta _k=\min\left(\frac{\left|
\lambda_{k+1}-\lambda _k\right| }2, \frac{\left|\lambda
_{k-1}-\lambda _k\right|}2\right)$.}}
\end{center}


\section{Introduction and notations.}
In this paper we deal with the spectral decomposition of a linear
operator $A$ over a Banach space $E$ (see for instance
\cite[Chapter 4]{ABHN} and \cite{AS}): For which subspace
$D\subset E$ do we have:
\begin{equation}\label{e1}
\sum_{k\in\mathbb{Z}}P_kx=x, \qquad (x\in D)
\end{equation}
where $\sigma(A):=\{i\lambda_k,k\in\mathbb{Z}\}$,
$\lambda_i\ne\lambda_j$, is assumed to be discrete and $P_k$ are
the associated eigen-projectors:
\begin{equation}\label{projector}
P_kx:=\int_{\gamma_k}(\lambda-A)^{-1}x\,d\lambda,
\end{equation}
whith $\gamma_k:=C(i\lambda_k,r_k)$ is the positively oriented
circle centered at $i\lambda_k$ and $r_k$ sufficiently small so
that $|\lambda_k-\lambda_{\ell}|>r_k$ for $\ell\ne k$.

\bigskip

In the case of bounded $C_0$-groups, if $\sup_k\|P_k\|<\infty$ and
if $\sum_k|\lambda_k|^{-1}<\infty$ then (\ref{e1}) holds true with
$D= D(A)$.

For the generator $A$ of a bounded $C_0$-group on a Banach space
$E$, there exists a large literature of spectral theory. The
introduction of the spectral theory was strongly motivated by
operator algebras (see for instance \cite{PO}), and was further
developed in the course of applications to groups of automorphisms
of operator algebras (see for instance \cite{BKR}). Always in the
framework of the generator of bounded groups, Gelfand's theorem
(1941) \cite{G} gives the equivalence between the spectrum of $A$
is zero, $\sigma(A)=\{0\}$
and $A$ is trivial, $A=0$.\\
Recently, related topics to Gelfand's theorem has been extensively
developed in a more general framework than bounded groups. (see
\cite{AR}, \cite[Chapter 4]{ABHN}, \cite{FJ}, \cite{Z} and the
references therein).

Unfortunately the framework of bounded groups fails for several
applications. A typical example is the Schr\"odinger operator
$i\Delta$ on $L^p(\Omega)$ with Dirichlet or Neumann boundary
condition for $p\ne 2$. In fact, for bounded domains
$\Omega\subset \mathbb{R}^N$ and using Sobolev embedding, Arendt
showed in 1991 \cite{Ar2} that $i\Delta$ (with Dirichlet or
Neumann boundary condition) generates a $k$-times integrated
group, that is polynomially bounded, for $k>\frac N2|\frac
1p-\frac 12|$. It should be noticed also that El-Mennaoui and
Kyantuo showed that the order $\frac N2|\frac 1p-\frac 12|$ is
optimal, see for instance \cite{EK}, and hence $i\Delta$ does
generate a $C_0$-group in this case.

The notion of integrated semigroups was introduced by Arendt
\cite{Ar} in 1987. An operator $A$ generates a $n$-times
integrated semigroup $\{S(t)\}$ ($n\in\mathbb{N}$) if and only if
there exists $w\in\mathbb{R}$ such that
$$(\lambda -A)^{-1}=\lambda^n \int_0^\infty e^{-\lambda t}S(t)\,dt$$
for $Re \lambda>w$. The $n$-times integrated group is said to be
polynomially bounded if there exists a positive integer $k$ and a
positive constant $C$ such that $\|S(t)\|\le C(1+|t|^k)$ for all
$t$. The $n$-times integrated group is said to be temperate if
there exists a positive constant $C$ such that $\|S(t)\|\le
C|t|^n$.

Except spectral distributions, there is no corresponding spectral
theory for integrated groups. The notion of spectral distributions
have been introduced in 1991. The special case of temperate
integrated groups has been extensively studied. See for instance
\cite{BEJ,J1,J2,EJ, FJ}.

In this paper, and using elementary calculations, we will show
that the decomposition (\ref{e1}) holds true for all $x\in
D(A^\ell)$ for some power $\ell$, where $A$ is the generator of
polynomially bounded $n$-times integrated group with discrete
spectrum.

Throughout this paper we suppose in addition (rescaling $A$ if
needed) that $(\lambda _k)_{k\in {\mathbb Z}^*}$ is an increasing
sequence, $\lambda _k>0$ for $k>0$ and $\lambda_k<0$ for $k<0$.

The paper is organized as follows: in the second section we prove
a separation result that permits to find two closed subspaces
$E_-$ and $E_+$ invariant by the integrated group such that
$\sigma(A_{|_{E_-}})\subset i\mathbb{R}^*_-$ and
$\sigma(A_{|_{E_+}})\subset i\mathbb{R}^*_+$. In the third section
the bounded group case is considered. The fourth section is
devoted to the integrated group case (Theorem \ref{t3}) with some
additional conditions on $|\lambda_k|^{-1}$ and on the distance
separating successive elements of the spectrum. In the last
section we apply our results to Schr\"odinger operators on $L^p$
spaces.

\textbf{Acknowledgment}. The third author wishes to thank A. El
Soufi and B. Helffer for enlighten discussions concerning this
paper.

\section{Some general separation results}

Given a linear operator $T$ on a Banach space $E$ with spectrum
$\sigma(T)=\sigma_1\cup \sigma_2$, union of two closed disjoint
sets with $\sigma_2$ compact, it is well known that
$P_2x:=\frac{1}{2i\pi} \int_\Gamma (\lambda-T)^{-1}x\,d\lambda$
$(x\in E)$, where $\Gamma$ is any bounded path surrounding
$\sigma_2$ and not containing in its interior any point of
$\sigma_1$, defines a bounded projector and then we have the
following decomposition: $E=E_1\oplus E_2$, where
$E_2:=P_2E=Im\,P_2$ and $E_1:=ker\,P_2$, and if we denote by
$T_1:=T_{|E_1}$ and $T_2:=T_{|E_2}$, the part of $T$ on $E_1$ and
$E_2$ respectively, then $\sigma(T_2)=\sigma_2$ and
$\sigma(T_1)=\sigma_1$. But if $\sigma_1$ and $\sigma_2$ are not
bounded, the associated projectors are not necessarily bounded.
Consider, for example, the operator $T=\frac d{dx}$ on the space
$E=L^1([0,2\pi])$. It is well known that the spectrum is
$\sigma(T)= i{\mathbb Z}$ and the projector defined by
$$P\left[ \sum_{k=-n}^n a_ke^{ikx}\right]:=\sum_{k=0}^n a_ke^{ikx}$$
is not bounded on $E$ (see for example \cite[Part II, Chapter 3
e.]{LT}).

The aim of this section is to show the following ``separation"
theorem:

\bigskip

\begin{theorem}\label{t1} Let $A$ be the generator of an $n$-times integrated
group $(S(t))_{t\in {\mathbb R}}$ satisfying $\|S(t)\|\leq
M(1+|t|^m)$ ($m\in \mathbb{R}_+$). Assume that $\sigma (A)=\left\{
i\lambda _n;\, n\in Z^*\right\} =\sigma^+\cup \sigma^-$ where
$\sigma^+:=\left\{ i\lambda_n, \,n\geq 1\right\}\subset i{\mathbb
R}_+^*$ and $\sigma^-:=\left\{ i\lambda_n,\,n\leq -1\right\}
\subset i{\mathbb R}_-^*$. Then there exists two closed subspaces
$E_-$ and $E_+$ of $E$ invariant by $(S(t))_{t\in {\mathbb R}}$
such that $D(A)\subset E_-\oplus E_+$ and $A_-:=A_{\mid_{E_-}}$
(resp. $A_+:=A_{\mid _{E_+}}$) generates $n$-times integrated
group satisfying the same boundedness property and whose spectrum
$\sigma (A_-)=\sigma ^-$ (resp. $\sigma (A_+)=\sigma ^+$).
\end{theorem}

\noindent For all $k\in {\mathbb Z}^*$, define the spectral
projection $P_k$ by
$$P_k:=\frac{1}{2i\pi}\int_{\gamma_k}R(\lambda, A)\,d\lambda,$$
where $\gamma_k:=\{z\in{\mathbb C};\quad
|z-i\lambda_k|=\varepsilon_k\}$ with $\varepsilon_k>0$ such that
$\overline{D}(i\lambda_k,\varepsilon_k)\cap
\sigma(A)=\{i\lambda_k\}$. It is clear that $P_k$ is a bounded
projection, $P_kx\in D(A)$ for all $x\in E$, and $P_kP_{k'}=0$ for
all $k\neq k'$.

\bigskip

We will need the following generalization of Gelfand's theorem,
for the proof we refer to \cite[Proposition 8.1]{El}.

\begin{proposition}\label{p1}\cite[Proposition 8.1]{El}\\
Let $A$ be the generator of a $n$-times integrated group
$(S(t))_{t\in{\mathbb R}}$ satisfying $\|S(t)\|\leq M(1+|t|^m)$,
with $m\geq n$ then $\sigma(A)=\emptyset$ if and only if
$E=\{0\}$.
\end{proposition}

\begin{remark}\label{r1} This implies the following weaker version
(see also \cite[Proposition 8.1]{El}): $\sigma(A)=\{0\}$ if and
only if $A^{m-n+1}=0$.
\end{remark}

In the following we are concerned with the convergence of the series
$\sum P_k x$. A first result is:

\begin{proposition}
$\displaystyle F:=\bigcup_{k\in {\mathbb Z}^*}P_k(E)$ is dense in $E$.
\end{proposition}
\noindent {\bf Proof:} Setting $G:=E/\overline{F}$, it suffices to
see that this Banach space is trivial. Denote by $B$ the part of
$A$ in $G$. $B$ generates an $n$-times integrated group that is
polynomially bounded as $(S(t))$ and the spectrum $\sigma
(B)=\emptyset$. Using proposition \ref{p1} $G=\{0\}$.
$\hfill\square$

\bigskip

Given $x\in F$, there exists $y\in E$, $k\in{\mathbb Z}^*$  such
that $x=P_k y$, then $P_\ell x=0$ for all $\ell\neq k$, thus we
obtain the convergence of the series $\sum P_k x$ to $x$ for all
$x\in F$.\\
In order to give more precise information on the
convergence of the series, we start by showing that the
convergence of the series $\sum_{\mathbb Z^*}P_kx$ is equivalent
to the convergence of $\sum_{k>0}P_kx$ and $\sum_{k<0}P_kx$ for
some regular $x$.

\bigskip

\begin{lemma}\label{l2}
For all $x\in E$ and all $k\in {\mathbb Z}^*$, the mapping $\lambda
\longmapsto R(\lambda ,A)P_kx$ is analytic on ${\mathbb C}\backslash\{i\lambda_k\} $
and there exists $c_k>0$ such that for all $\lambda \notin \overline{D}
(i\lambda _k,\varepsilon _k)$
$$
\left\| R(\lambda ,A)P_k\right\| \leq \frac{C_k}{dist(\lambda ,\gamma _k)}.
$$
\end{lemma}
\noindent {\bf Proof :} Let $k\in {\mathbb Z}^*$ and $x\in E$.
Since $R(\lambda ,A)P_kx=\frac 1{2i\pi} \int_{\gamma _k}R(\lambda
,A)R(\mu ,A)x\,d\mu$, the resolvent equation implies
\begin{eqnarray*}
R(\lambda ,A)P_kx&=&\frac 1{2i\pi }\int_{\gamma _k}\frac{R(\lambda ,A)x-R(\mu
,A)x}{\mu -\lambda }\,d\mu\\
&=&\frac 1{2i\pi }\int_{\gamma _k}\frac{R(\lambda ,A)x}{\mu
-\lambda }\,d\mu -\frac 1{2i\pi }\int_{\gamma _k}\frac{R(\mu ,A)x}{\mu
-\lambda }\,d\mu .
\end{eqnarray*}
\noindent For $\lambda \notin \overline{D}(i\lambda _k,\varepsilon _k)$,
$\displaystyle \frac 1{2i\pi }\int_{\gamma _k}\frac{R(\lambda ,A)x}{\mu -\lambda }
\,d\mu =0$, then
$$
R(\lambda ,A)P_kx=-\frac 1{2i\pi }\int_{\gamma _k}\frac{R(\mu ,A)x}{\mu
-\lambda }\,d\mu,
$$
thus $\displaystyle\left\| R(\lambda ,A)P_kx\right\| \leq
\frac{\left\| x\right\| \varepsilon _k}{dist(\lambda ,\gamma
_k)}\sup_{\mu \in \gamma _k}\left\| R(\mu ,A)\right\|$. It
suffices then to take $\displaystyle C_k:=\varepsilon _k\sup_{\mu
\in \gamma _k}\left\| R(\mu ,A)\right\|$. $\hfill\square$

\begin{proposition}
There exists ${\mathbb P}\in \mathcal{L}(D(A),E)$ such that
\begin{equation}\label{eq3}
{\mathbb P}\left[\sum_{-N\leq k\neq 0\leq
N}A^{-m-1}P_{k}x\right]=\sum_{k=1}^{k=N}A^{-m-1}P_{k}x \qquad
(N\in {\mathbb N}^*,\,x\in E),
\end{equation}
$$
{\left({\mathbb P}A^{-m-1}\right)}^2={\mathbb P}A^{-2m-2}.
$$
\end{proposition}
\noindent {\bf Proof:} Let $\delta >0$ be such that $\delta <\min
(\lambda _{1},\left| \lambda _{-1}\right| )$ and consider
$(R_n)_{n\in {\mathbb N}^*}$ an increasing sequence of strictly
positive numbers satisfying for all $N\geq 1$
$$
\overline{D}(0,R_N)\cap \sigma (A)=\left\{ i\lambda_k;\qquad -N\leq k\leq N,\quad k\ne0\right\} .
$$
Let $\Gamma _N$ be the positively oriented path
$\Gamma_N:=[-R_N,-\delta ]\cup C^+(0,\delta )\cup [\delta
,R_N]$ (see the figure below).\\
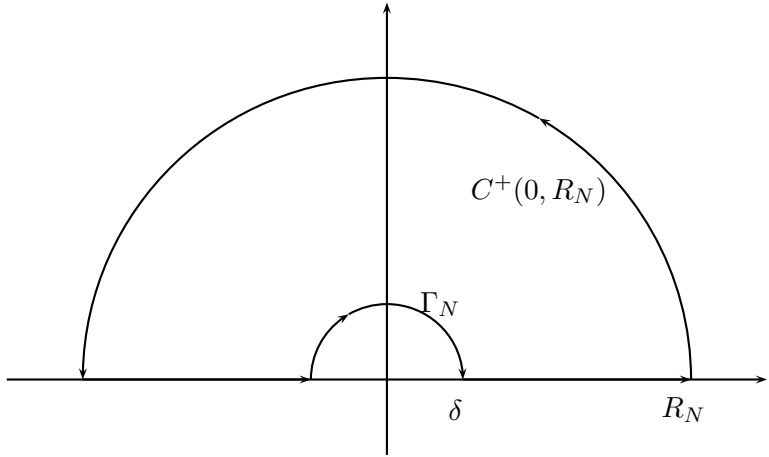
\begin{center}
%figure path
\pspicture(-5,-1.5)(5,5.5) \psline{->}(-5,0)(5,0)
\psline{->}(0,-1)(0,5) \psline{->}(1,0)(4,0)
\psline{->}(-4,0)(-1,0) \psarc{<-}(0,0){1}{0}{120}
\psarc{<-}(0,0){1}{120}{180} \rput(.9,-.4){$\delta$}
\rput(.7,1){$\Gamma_N$} \rput(3.9,-.4){$R_N$}
\psarc{->}(0,0){4}{0}{60} \psarc{->}(0,0){4}{60}{180}
\rput(2,2.5){$C^+(0,R_N)$}
\endpspicture
\end{center}

\noindent Since $0\not \in \sigma(A)$, consider, for $N\in{\mathbb
N}^*$ and $x\in E$, the bounded operator $\displaystyle
Q_NA^{-m-1}x:=\frac 1{2i\pi} \int_{\Gamma _N}R(\lambda
,A)A^{-m-1}x\,d\lambda$. We will first show that $Q_NA^{-m-1}$
converges strongly in $E$. The resolvent equation:
$$R(\lambda,A)A^{-k}x=\frac{A^{-k}x}\lambda+\cdots+\frac{A^{-1}x}{\lambda^k}+
\frac{R(\lambda,A)x}{\lambda^k}$$
implies
\begin{eqnarray*}
Q_NA^{-m-1}x=\frac{A^{-m-1}x}2&+&\sum_{j=2}^{m+1}\frac 1{2i\pi }\int_{\Gamma _N}
(-\lambda)^{-j}A^{-(m+1-j)}x\, d\lambda\\ &+& \frac 1{2i\pi }\int_{\Gamma _N}
(-\lambda)^{-m-1}R(\lambda ,A)x\, d\lambda .
\end{eqnarray*}
Since $\int_{\Gamma _N} (-\lambda)^{-j}A^{-(m+1-j)}x\, d\lambda = -\int_{C^+(0,R_N)}
(-\lambda)^{-j}A^{-(m+1-j)}x\, d\lambda$, whose modulus tends to zero as $N\to\infty$
for $j>1$, it suffices to prove the convergence of the last term. Using the bound
on the integrated group $(S(t))$ we have $(Re\lambda\neq 0)$
$$\left\|\frac{R(\lambda, A)}{\lambda^{m+1}} \right\|\leq C|\lambda|^{n-m+1}
\left(\frac 1{|Re \lambda|}+\frac 1{|Re \lambda|^{m+1}}\right)$$
so $\displaystyle \lim_{N\to\infty}\frac 1{2i\pi }\int_{\Gamma _N}
R(\lambda ,A)x(-\lambda)^{-m-1}\, d\lambda$ exists. Denote by
$QA^{-m-1}x$ its limit. Set ${\mathbb P}:=Q-\frac 12Id$ and
let's show the equality (\ref{eq3}).\\
From one side, for $1\leq k\leq N$, by Cauchy formula and Lemma \ref{l2} we have:
\begin{eqnarray*}
Q_NA^{-m-1}P_kx&=&\frac 1{2i\pi}
\int_{\Gamma _N}R(\lambda ,A)A^{-m-1}P_kx\,d\lambda\\
&=&\int_{C^+(0,R_N)}R(\lambda ,A)A^{-m-1}P_kx\,d\lambda+P_kA^{-m-1}P_kx\\
&=&-\frac 1{2i\pi }\int_{C^+(0,R_N)}
\sum_{j=1}^{m+1} (-\lambda)^{-j}A^{-(m+1-j)}P_kx\,d\lambda\\
&&-\frac 1{2i\pi }\int_{C^+(0,R_N)}(-\lambda)^{-m-1} R(\lambda
,A)P_kx\, d\lambda+A^{-m-1}P_kx\\
&=&A^{-m-1}P_kx+\frac 1{2i\pi
}\int_{C^+(0,R_N)}\frac{R(\lambda
,A)P_kx+A^{-1}P_kx}\lambda\, d\lambda \\
&\longrightarrow&A^{-m-1}P_kx+\frac{A^{-m-1}P_kx}2\quad\quad (\mbox{as } N\to\infty)\\
&{}&=\frac 32A^{-m-1}P_kx=QA^{-m-1}P_kx.
\end{eqnarray*}

\noindent Similarly, for $-N\leq k\leq -1$, we get\\
$$\lim_{N\rightarrow +\infty}Q_NA^{-m-1}P_kx=QA^{-m-1}P_kx =\frac
12A^{-m-1}P_kx.$$ Thus we have proved the equality (\ref{eq3}).
Now since $F$ is dense in $E$ and ${\mathbb P}A^{-1}\in {\cal
L}(E)$, we see that $\displaystyle{\left({\mathbb
P}A^{-m-1}\right)}^2 ={\mathbb P}A^{-2m-2}$. $\hfill\square$

\bigskip

\noindent{\bf Remarks 2}\\
1- Notice that $\displaystyle {\mathbb P}x=-\frac 1{2i\pi}\int_\Gamma
\frac{R(\lambda ,A)A^{m+1}x}{(-\lambda)^{m+1}}\, d\lambda$ for all $x\in D(A^{m+1})$,
where
$$
\Gamma:=]-\infty ,-\delta]\cup C^+(0,\delta )\cup [\delta ,+\infty [.
$$

\noindent 2- In general the operator ${\mathbb P}$ is not bounded.

\noindent For example consider the rotation group
$(T(t))_{t\in{\mathbb R}}$ on $L_{2\pi }^p({\mathbb R})$ with
$1\leq p<+\infty $, where $L_{2\pi }^p({\mathbb R})$ the Banach
space of measurable $2\pi$-periodic functions on ${\mathbb R}$
satisfying $\int_0^{2\pi}|f(x)|^p\,dx<+\infty $. The rotation
group is defined by $T(t)f(x):=f(x+t)$ for $f\in
L_{2\pi}^p({\mathbb R})$. The generator $A$ of
$(T(t))_{t\in{\mathbb R}}$ is defined by $Af=f'$ and $D(A)=\{ f\in
L_{2\pi }^p({\mathbb R})$,  $f'\in L_{2\pi }^p({\mathbb R}) \}$.
$f'$ is the derivative in the distribution sense. It is easy to
verify that the spectrum of $A$, $\sigma (A)=i{\mathbb Z}$ and the
eigenvectors are exactly the Fourier basis $(e^{int})_{n\in
{\mathbb Z}}$. Therefore ${\mathbb P}{(f)}$ is the Fourier series
associated to $f$. It is well known that for $p\neq 2$,
$(e^{int})_{n\in {\mathbb Z}}$ is not a basis of $L_{2\pi
}^p({\mathbb R})$ so ${\mathbb P}{(f)}$ does not converge in
$L_{2\pi }^{p}({\mathbb R})$ except for $p=2$. Consequently
${\mathbb P}$ is defined only on the domain of $A$.

\bigskip

\noindent {\bf Proof of Theorem \ref{t1}:} Set $E_-:=\{ x\in E;\,
{\mathbb P}A^{-m-1}x=0\} = ker {\mathbb P}A^{-m-1}$ and $E_+:=\{
x\in E;$ ${\mathbb P}A^{-m-1}x=A^{-1}x\}$. $E_-$ and $E_+$ are two
closed subspaces of $E$. It is clear that $D(A^{m+1})\subset
E_0+E_1$ and that $E_0$ and $E_1$ are invariant by $(S(t))_{t\in
{\mathbb R}}$. If $x\in E_0\cap E_1$, then $A^{-m-1}x={\mathbb
P}A^{-m-1}x=0$ and so $x=0$.

\noindent Let $A_-$ be the part of $A$ in $E_-$. It's clear that
$\sigma (A_-)\subset \sigma (A)$. To show that $\sigma
(A_-)=\sigma ^-$ it suffices to verify that $\sigma (A_-)\cap
\sigma ^+=\emptyset$. Suppose that there exists $k\in {\mathbb
N}^*$ such that $i\lambda _{k}\in \sigma(A_-)$. There exists then
a sequence $(x_{n})_{n\in {\mathbb N}^*}\subset D(A_{-})$ such
that for all $j\in {\mathbb N}^*$ we have $\| x_j\| =1$ and
$\left\| Ax_j-i\lambda _kx_j\right\| \leq \frac 1j$.

\noindent From one side for $\lambda \in \Gamma$ we have
$$
\left\| \frac{-x_j}{\lambda -i\lambda _k}+R(\lambda ,A)x_j\right\| =
\left\| \frac{R(\lambda ,A)(Ax_j-i\lambda _kx_j)}{\lambda -i\lambda _k}
\right\| \leq \frac C{j\left| \lambda -\lambda _k\right| }
$$
where $c:=\sup_{\lambda \in \Gamma }\left\| R(\lambda ,A)\right\|$.\\ From the other side
we have $\displaystyle{\mathbb P}x=\frac{1}{2i\pi }\int_{\Gamma }\frac{R(\lambda ,A)
A^{m+1}x}{(-\lambda)^{m+1} }\,d\lambda $ for $x\in D(A^{m+1})$ and since
$\frac 1{2i\pi}\int_{\Gamma
}\frac{x_j}{(-\lambda)^{m+1} (\lambda -i\lambda _k)}\,d\lambda =\frac{x_j}
{(-i\lambda _k)^{m+1}}$, then
$$
\left\|{\mathbb P}A_-^{-m-1}x_j-\frac{x_j}{(-i\lambda _k)^{m+1}}\right\|
\leq \frac 1{2\pi }\left|\int_\Gamma \frac C{j \lambda (\lambda
-i\lambda _k)}\,d\lambda\right|.$$
Thus $\displaystyle \lim_{j\to\infty} \left\|{\mathbb P}A_-^{-m-1}x_j-
\frac{x_j}{(-i\lambda _k)^{m+1}}\right\|=0$, but $x_j\in E_-$, so
${\mathbb P}A^{-m-1}x_j=0$ hence $\lim_{n\to +\infty}\| x_j\| =0$. This is absurd since
$\left\| x_j\right\| =1$. Therefore $\sigma (A_-)\cap \sigma ^+=\emptyset$ and the
spectrum of $A_-$ is $\sigma ^-$. $\hfill\square$

\section{Bounded $C_0$-groups case}

The case of bounded groups is more classical and spectral theory
must lead to simple proofs. Indeed, the spectral projections can
be written in the form
$P_k=\int_{\mathbb{R}}\mathcal{F}(f_k)(t)T(t)\,dt$ where
$f_k(\lambda_j)=\delta_{kj}$. Assume for the moment that
$c:=\sup_k\|P_k\|<\infty$. For $x\in D(A)$, and using Gelfand's
theorem, one has $P_kAx=\lambda_k P_kx$ and hence $\|P_kx\|\le
c|\lambda_k|^{-1}\|Ax\|$. Under the hypothesis that
$\sum_k|\lambda_k|^{-1}<\infty$ it follows that $\sum_kP_kx$ is
convergent for all $x\in D(A)$.\\
The problem now is to choose suitable functions $f_k$ so that
$\sup_k\|\mathcal{F}f_k\|_{L^1}$ is finite. For this take
$$f_k(t):=f\left(\frac{t-\lambda_k}{\delta_k}\right)$$ where $f$ is
any $C^\infty$ function of support in $(-1,1)$ that is equal to 1
on a neighborhood of zero, and $\delta_k:=\frac
12\min(\lambda_k-\lambda_{k-1},\lambda_{k+1}-\lambda_k)$. Then we
get
$$\|P_k\|\le \|\mathcal{F}f_k\|_{L^1}=\|\mathcal{F}f\|_{L^1}.$$
We have thus proved the following theorem. However, aiming to
introduce our method, we give another proof based on the contour
method.

\begin{theorem}\label{t2}
Let $A$ be the generator of a bounded $C_0$-group
$(T(t))_{t\in {\mathbb R}}$ with discrete spectrum
$\sigma(A)=\{i\lambda_k,\,k\in\mathbb{Z}^*\}$. Assume that $\sum
_{k\in {\mathbb Z}^*} \frac 1{\left|\lambda _k\right|}<\infty$,
then $\displaystyle\sum_{k\in {\mathbb Z}^*}P_kx=x$ $(x\in D(A))$
and $\displaystyle\sum_{k\in {\mathbb Z}^*}i\lambda _kP_kx=Ax$
$(x\in D(A^{2}))$.
\end{theorem}
\noindent {\bf Proof:} Let's show first that
$\sup_k\|P_k\|<\infty$, where $P_k$ is defined by
(\ref{projector}). Remark that we can write
$$P_kx=\frac 1{2i\pi}\int_{\gamma_k}(\lambda-A)^{-1}\left[
1+\frac{(\lambda-i\lambda_k)^2}{r_k^2}\right]\,d\lambda$$ since
$\lambda_k$ is a simple pole for $(\lambda-A)^{-1}$. Now since $A$
generates a bounded group we have $\|Re(\lambda)(\lambda-A)^{-1}\|
\le c$. Thus
\begin{eqnarray*}
\|P_k\|&\le&\frac 1{2\pi}\int_0^{2\pi}\|(r_ke^{i\theta}
-A)^{-1}\|\, |1+e^{2i\theta}|\,r_kd\theta\\
&\le&\frac 1{\pi}\int_0^{2\pi}\|r_k\cos\theta(r_ke^{i\theta}
-A)^{-1}\|\,d\theta\le 2c.
\end{eqnarray*}
Now, for $x\in D(A)$ and using Gelfand's theorem, we have
$$P_kx=\frac1{\lambda_k}AP_kx=\frac1{\lambda_k}P_kAx$$
hence
$$\|P_kx\|=\frac1{|\lambda_k|}\|P_kAx\|\le \frac{2c}{|\lambda_k|}\|Ax\|.$$
Therefore, $\sum_kP_kx$ is normally convergent for all $x\in
D(A)$. Since $\sum_kP_kx=x$ for all $x\in F$ we get the result
using the density of $F$ in $D(A)$. $\hfill\square$

\section{Integrated groups case}

As we saw in the preceding section, the spectral theory could give
easily the spectral decomposition in the case of bounded groups.
But this method does not work\footnote{Since the uniform
boundedness of the projectors is not true in general in this
case.} for the more general class, while contour method still
applicable.

In this section we will give more appropriate conditions on the
behavior of the spectrum to get convergence.

\bigskip
The following is an ``integrated" version of theorem \ref{t2}.
Here $n$ and $\ell$ are fixed integers.

\begin{theorem}\label{t3} Let $A$ be the generator of an $n$-times
integrated group $(S(t))_{t\in{\mathbb R}}$ satisfying for all $t$
\begin{equation}\label{pol}\|S(t)\|\le M(1+|t|^n),\end{equation}
with discrete spectrum
$\sigma(A)=\{i\lambda_k;\,k\in\mathbb{Z}^*\}\subset
i\mathbb{R}^*$. Assume that
\begin{equation}\label{somme}\sum
_{k\in {\mathbb Z}^*} \frac 1 {\left|\lambda _k\right|^\ell
\delta_k^n}<\infty,\end{equation} where
$$\displaystyle\delta_k:=\min \left(\frac{| \lambda _{k+1}-\lambda
_k|}2,\frac{| \lambda _{k-1}-\lambda _k|}2 \right),$$ then
\begin{equation}\label{decomposition}
\sum_{k\in {\mathbb Z}^*}P_kx=x\qquad\mbox{for all }x\in
D(A^{n+\ell}).\end{equation}
\end{theorem}
\noindent {\bf Proof:} Using theorem \ref{t1}, we can suppose that
$\sigma (A)=\left\{ i\lambda _n;\quad n\in {\mathbb N}^*\right\}
\subset i {\mathbb R}_+^*$. Let $k\in {\mathbb Z}^*$ and set
$A_k:=(A-i\lambda _kId)\mid _{P_k(E)}$ the part of $A-i\lambda
_{k}Id$ on $P_k(E)$. $A_k$ is a bounded operator on $P_k(E)$,
hence it generates, in particular, an $n$-times integrated group
$(S_k(t))_{t\in {\mathbb R}}$ that is temperate at infinity (i.e.
$\|S_k(t\|=O(|t|^{2n}))$) and whose spectrum $\sigma (A_{k})=\{
0\}$. Hence, using Remark \ref{r1}, we have $A_k^{n+1}=0$. Thus we
have
$$
R(i\lambda _k+\lambda,A
)P_{\ell}x=\sum_{j=1}^n\frac{A_k^j}{\left(\lambda -i(\lambda
_\ell-\lambda _k)\right)^{j+1}}P_\ell x
$$
for all $\left| \lambda \right| =\delta_k$.

\noindent Observing that the mapping $\displaystyle \lambda
\longmapsto \lambda ^{2(2n+1)}R(i\lambda _k+\lambda )x$ is
analytic from the disc $\overline{D}(0, \delta_k)$ into $P_k(E)$,
denote by
$${\mathbb P}_Nx:=\frac 1{2i\pi }\int_{\left| \lambda \right|
=\delta_k }{\left[1+\frac{ \lambda ^{2(2n+1)}}{\delta_k
^{2(2n+1)}} \right]}^{n+1} \sum_{k=1}^{k=N}R(i\lambda _k+\lambda,
A )x\,d\lambda,
$$
then
$$
\left\{
\begin{array}{lr}
{\mathbb P}_NP_\ell(x)=P_\ell(x)& \mbox{for } |\ell| \leq N\\
{\mathbb P}_NP_\ell(x)=0& \mbox{for } |\ell| >N
\end{array}
\right.
$$
We deduce that for all $\displaystyle x\in F:=\bigcup_{k>1}P_k(E)$
the sequence $({\mathbb P}_Nx)_{N\in {\mathbb N}^*}$ converges to $x$.

\noindent Using the following resolvent equation:
$$R(\lambda,A)A^{-n}x=\frac{A^{-n}x}\lambda+\cdots+\frac{A^{-1}x}{\lambda^n}+
\frac{R(\lambda,A)x}{\lambda^n}$$
we get
$${\mathbb P}_NA^{-n}x=\frac 1{2i\pi }\int_{|\lambda| =\delta }{\left[1+
\frac{\lambda ^{2(2n+1)}}{\delta_k^{2(2n+1)}}\right]}^{n+1}\,
\sum_{k=1}^{k=N}\frac{ R(i\lambda _k+\lambda,A )}{(\lambda
+i\lambda _k)^n}x\,d\lambda
$$
Now for $Re(\lambda )>0$, we have $(\lambda +i\lambda
_k)^{-n}R(i\lambda _k+\lambda,A )x =\int_{0}^{+\infty}e^{-i\lambda
_kt}e^{-\lambda t}S(t)x\,dt$. Thus, for $x\in D(A^{n+\ell})$,
setting $U_N(\lambda)x:=\sum_{k=1}^{k=N} (\lambda +i\lambda
_k)^{-n}R(i\lambda _{k}+\lambda,A )x$, we have
\begin{eqnarray*}
U_N(\lambda )x&=&\sum_{k=1}^{k=N}\int_0^{+\infty }e^{-i\lambda
_kt}e^{-\lambda t}
S(t)x\, dt\nonumber\\
&=&\int_0^{+\infty
}\left[\sum_{k=1}^{k=N}e^{-i\lambda_kt}\right]e^{-\lambda
t}S(t)x\,dt
\end{eqnarray*}
Integrating $\ell$ times by parts, we get
\begin{eqnarray*}
U_N(\lambda )x&=&\sum_{k=1}^{k=N}\int_0^{+\infty }e^{-i\lambda
_kt}e^{-\lambda t}
S(t)x\, dt\\
&=&\int_0^{+\infty }\left[\sum_{k=1}^{k=N}e^{-i\lambda_kt}\right]e^{-\lambda t}S(t)x\,dt\\
&=&\int_0^{+\infty }\sum_{k=1}^{k=N}\frac{e^{-i\lambda
_kt}}{(i\lambda _k)^\ell} e^{-\lambda t}\left[S(t)(A-\lambda)^\ell
+f(\lambda,A,t)\right]x\,dt,
\end{eqnarray*}
where $f(\lambda,A,t):=\sum_{0\le\alpha,\beta\le
\ell}\sum_{0\le\gamma\le n} a_{\alpha,\beta,\gamma}
\lambda^{\alpha}t^{\gamma}A^{\beta}$ is a polynomial on $t$,
$a_{\alpha,\beta,\gamma}$ are constants. Using (\ref{pol}) we get
\begin{eqnarray*}
\| U_N(\lambda )x\| &\leq& c_1\| x\|_{A^\ell}
\sum_{k=1}^{k=N}\frac1{|\lambda_k|^\ell }\int_0^{+\infty}
(\sum_{0\le j\le n}t^j)e^{-t Re(\lambda )}\,dt\\ &\leq& c_2\|
x\|_{A^\ell} \left[\sum_{k=1}^{k=N}\frac1{|\lambda_k|^\ell
}\right] \sum_{1\le j\le n+1}\frac 1{| Re(\lambda )|^j},
\end{eqnarray*}
where $c_1$ and $c_2$ are two strictly positive constants. Then,
denoting by $C_k^+=\{Rez>0,\,|z|=\delta_k\}$ and by ${\mathbb
P}_N^+$
$${\mathbb P}_N^+x:=\frac 1{2i\pi }\int_{C_k^+}{\left[1+\frac{
\lambda ^{2(2n+1)}}{\delta_k^{2(2n+1)}}\right]}^{n+1}
\sum_{k=1}^{k=N}\frac{R(i\lambda _k+\lambda,A
)}{(\lambda+i\lambda_k)^n}x\,d\lambda,
$$
we get
\begin{eqnarray*}
\left\| {\mathbb P}_N^+A^{-n-1}x\right\|&\le& c_3\delta_k
\|x\|_{A^\ell}\sum_{k=1}^{k=N}\frac1{|\lambda_k|^\ell}\int_{-\pi/2}^{\pi/2}
\left| 1+e^{2i\theta (2n+1)}\right|^{n+1}\times\\ &&{} \hskip 4cm
\sum_{1\le j\le n+1} \frac 1{\left| \delta_k \cos \theta \right|
^{j}}\,d\theta\\ &\leq& c_5\|x\|_{A^\ell}\sum_{0\le k\le
N}\sum_{0\le j\le n}\frac1{|\lambda_k|^\ell \delta_k^j}
\int_{-\pi/2}^{\pi/2} \frac {|\cos(\theta (2n+1))|^{n+1}} {|\cos
\theta|^{j+1}}\,d\theta.
\end{eqnarray*}
Thus (we may assume that $\delta_k\le 1$)
\begin{eqnarray*}
\left\| {\mathbb P}_NA^{-n-1}x\right\|\le
c_5\|x\|_{A^\ell}\sum_{0\le k\le N}\frac1{|\lambda_k|^\ell
\delta_k^n} \int_0^{2\pi}\sum_{0\le j\le n} \frac {|\cos(\theta
(2n+1))|^{n+1}} {|\cos \theta|^j}\,d\theta.
\end{eqnarray*}
Since the integrals $$\displaystyle \int_0^{2\pi }\frac{\left|
\cos \theta (2n+1) \right| ^{n+1}}{\left| \cos\theta \right|
^{j}}\,d\theta = 2\int_{-\frac{\pi }{2}}^{+\frac{\pi
}{2}}\frac{\left| \sin \theta (2n+1)\right| ^{n+1}} {\left| \sin
\theta \right|^{j}}\,d\theta<\infty $$ are finite for $j$, and
using (\ref{somme}) there exists $c>0$ such that for all $x\in E$
we have
$$
\left\| {\mathbb P}_NA^{-(n+1)}x\right\| \leq c\left\| x\right\|.
$$
But for $x\in F$, $ \lim_{N\rightarrow +\infty }{\mathbb P}_Nx=x$
and $F$ is dense in $E$ then for $x\in D(A^{n+1})$
$\lim_{N\rightarrow +\infty } {\mathbb P}_Nx=x$. $\hfill\square$

\section{Applications: Schr\"odinger operator on a compact manifold}

Consider the Sch\"odinger operator $H_p=-\Delta+V$, where $\Delta$
is the Laplace-Beltrami operator on $L^p(\mathcal{M})$,
$\mathcal{M}$ is a compact manifold with Dirichlet or Neumann
manifold in the nonempty boundary case and $V$ a regular
potential. The operator $H_2$ generates on $L^2(\mathcal{M})$ a
self-adjoint operator. Following the same technique as in
\cite[Theorem 4.3, p.37]{Ar2} and by Sobolev imbeddings, the
operator $iH_p$ generates on $L^p(\mathcal{M})$ a $n$-times
integrated group $\{S(t)\}$, for $n> \frac N2|\frac 1p-\frac 12|$,
with $1<p<\infty$. Now by a similar calculation as in
\cite[Proposition 3.1, p. 61]{J3} the integrated group is
polynomially bounded: $\|S(t)\|\le C(1+|t|^n)$ for some positive
constant $C$ and all $t\in\mathbb{R}$.

The basic fact in this section is the $p$-spectral independence of
the Schr\"odinger operator on $L^p$ for $1\le p<\infty$ (see
\cite{HV}, \cite{Ar3} and the references therein). Moreover, it's
well known that the spectrum of $-\Delta+V$ is a strictly
increasing sequence $(\lambda_k)$ and if $(\mu_k)$ are  the
eigenvalues counted with their multiplicities then $\mu_k\sim
Ck^{2/N}$ as $k\to\infty$ where $C$ is a positive constant
independent of $k$. Indeed, for any $k$ we have
$$\mu_k(-\Delta)-\|V\|_\infty\le \mu_k(-\Delta+V)\le
\mu_k(-\Delta)+\|V\|_\infty.$$

To apply Theorem \ref{t3} we need an estimate for
$\lambda_{k+1}-\lambda_k$.

In what follows we give applications for which one can use our
result. Unfortunately, in many cases the condition (\ref{somme})
is not satisfied (which shows the relevance of this condition).

\subsection{Schr\"odinger operator on $\mathbf{L^p(S^N)}$}

Consider the Laplace-Beltrami operator $-\Delta$ on the sphere
$S^N$. It is well known (see for instance \cite[Corollary
4.3]{Tay}) that the spectrum of the Laplace operator $-\Delta$ on
the space $L^2(S^N)$ (and hence on $L^p$) is the set $\{k(k+N-1),
\, k\in \mathbb{Z}\}$. In this case one can apply Theorem \ref{t3}
with $\ell>1/2$ and we have the decomposition
(\ref{decomposition}).

It is well known that the dimension of the eigen-space associated
to the eigenvalue $k(k+N-1)$ is $m(k):=(n+2k-1)\frac{(N+k-2)!}
{k!(N-1)!}$. Each eigenvalue $k(k+N-1)$ of $-\Delta$ is split into
$m(k)$ eigenvalues of the operator $H_p:=-\Delta+V$, with some
regular potential $V$ and the distance between two successive
eigenvalues is now the distance between the two successive
$\mu_{k,\ell}$ and $\mu_{k,\ell+1}$ (supposed ordered). Several
authors studied the asymptotic of the $(\mu_{k,\ell})_{0\le
\ell\le 2k}$, and in particular, Grigis gave in the two
dimensional case (see \cite{Grigis}) a very simple form of
potential $V(x)=4x_1x_2$, with the notation
$S^2=\{x_1^2+x_2^2+x_3^2=1\}$, for which the difference
$\mu_{k,\ell+1}-\mu_{k,\ell}=O(k^{-\infty})$ is exponentially
small. See also \cite[Theorem 7.9]{CV} for another example. This
shows the relevance of the hypothesis (\ref{somme}).

\subsection{Harmonic oscillator on $\mathbf{L^p(\mathbb{R}^N)}$}

Consider the operator $H_p=-\Delta+|x|^2$ on $L^p(\mathbb{R}^N)$.
For $p=2$, $H_2$ is self-adjoint and its spectrum is $\{ 2k+N;\,
k\in\mathbb{N}\}$ (see for instance \cite[Chapter 8, section
6]{Tay}). In this case the distances between successive
eigenvalues are bounded below. Using the same technics as in the
last subsection Theorem \ref{t3} is applicable and we have the
decomposition
$$\sum_{k\ge 1}P_kx=x\qquad (x\in D(H^{n+\ell})),$$
where $n>\frac N2|\frac 1p-\frac 12|$ and $\ell>1$.

\subsection{Schr\"odinger operator on a 2-D flat torus}

Consider the Laplace-Beltrami operator $H=-\Delta_{{T}}$ on the
torus ${T}:=\mathbb{R}^2/\Gamma$ where
$\Gamma:=\mathbb{Z}e_1\oplus\mathbb{Z}e_2$ and $e_1=(a,0)$,
$e_2=(0,b)$, $a$ and $b$ are positive real numbers. Using again
Sobolev imbeddings, the operator $iH$ generates on $L^p({T})$ a
$n$-times integrated group $\{S(t)\}$, for $n\ge \frac N2|\frac
1p-\frac 12|$, with $1<p<\infty$. On $L^p({T})$, $1\le p<\infty$,
the spectrum $\sigma(H)$ is formed of eigenvalues and is given by
$\{a'm^2+b'n^2;\,m,n\in\mathbb{Z} \}$ (see for instance
\cite[Chapitre 3]{BGM}), where $a':=(2\pi/a)^2$ and
$b':=(2\pi/b)^2$. In order to apply Theorem \ref{t3} we need to
estimate the decay of $\delta_k$, the distance between successive
eigenvalues. For this we need the following algebraic proposition:

\begin{proposition}\cite[Proposition D.1.2 and Remark D.1.2.2]{HS}
Let $\alpha$ be an algebraic number of degree $d\ge 2$. There
exists a constant $C(\alpha)>0$ such that $$\left|\alpha -\frac
pq\right|\ge \frac {C(\alpha)}{q^d}\hskip 1cm \mbox{for all }\,
\frac pq\in\mathbb{Q}$$
\end{proposition}

\begin{corollary}
Let $\alpha:=b^2/a^2$ and $d$ be as in the last proposition and
denote by $\{\lambda_n,\,n\in\mathbb{N}\}$ the spectrum of
$i\Delta_T$ the Schr\"odinger operator on $L^p(T)$, $1\le
p<\infty$. Then there exists projectors $P_k$ such that
$$\sum_{k\ge 1}P_kx=x\qquad (x\in D((-\Delta)^{n(d-1)+2})),$$
where $n>\frac N2|\frac 1p-\frac 12|$.
\end{corollary}
\noindent\textbf{Proof.} Denote by $0<\bar \lambda=a'\bar
m^2+b'\bar n^2<\lambda=a'm^2+b'n^2$ be two successive large
eigenvalues. Then $$[\lambda-\bar \lambda]/b'=|n^2-\bar
n^2+\alpha(m^2-\bar m^2)|=|m^2-\bar m^2|\,\left|\frac{n^2-\bar
n^2}{m^2-\bar m^2}+\alpha\right|$$ and by the last proposition
\begin{eqnarray*}
\frac{\lambda-\bar \lambda}{b'}\ge \frac {C(\alpha)}{|m^2-\bar
m^2|^{d-1}}\ge \frac {C}{\lambda^{d-1}}.\end{eqnarray*} Therefore
$\lambda_k^{d-1}\delta_k\ge const$ (with the notations of Theorem
\ref{t3}). Applying Theorem \ref{t3}, and since the series
$\sum_{n\ge 1}\sum_{m\ge 1}\frac 1{(am^2+bn^2)^2}$ is convergent,
one gets the desired result. $\hfill\square$

%============================================================================

\noindent A. Driouich and O. El-Mennaoui\\
Agadir university, Mathematics department\\
Agadri, Morocco\\
Email:driouichabder@hotmail.com\hskip 1cm omar1172@caramail.com\\

\bigskip

\noindent M. Jazar\\
Lebanese university, Mathematics department\\
P.O. box 155-012, Beirut Lebanon\\
Email: mjazar@ul.edu.lb

\end{document}